\documentclass[12pt]{article}
\usepackage{fullpage}
\usepackage{amssymb,amsmath,amsfonts,amsthm,latexsym,enumerate,url,cases}
\numberwithin{equation}{section}
\usepackage{hyperref}
\hypersetup{colorlinks=true,citecolor=blue,linkcolor=blue,urlcolor=blue}

\newtheorem{theorem}{Theorem}[section] %
\newtheorem{problem}[theorem]{Problem} %

\begin{document}
\title{On $AP_3$ - covering sequences\footnote{The work is supported by the
National Natural Science Foundation of China, Grant No. 11771211
and a project funded by the Priority Academic Program Development
of Jiangsu Higher Education Institutions.}}

\author{  Yong-Gao Chen\footnote{E-mail address: ygchen@njnu.edu.cn}\\
\small  School of Mathematical Sciences and Institute of Mathematics, \\
\small  Nanjing Normal University,  Nanjing  210023,  P. R. China
}
\date{}
\maketitle \baselineskip 18pt \maketitle \baselineskip 18pt

{\bf Abstract.}  Recently, motivated by Stanley sequences, Kiss,
S\' andor and Yang introduced a new type sequence: a sequence $A$
of nonnegative integers is called an $AP_k$ - covering sequence if
there exists an integer $n_0$ such that if $n
> n_0$, then there exist $a_1\in A, \dots , a_{k-1}\in A$,
$a_1<a_2<\cdots <a_{k-1}<n$ such that $a_1, \dots , a_{k-1}, n$
form a $k$-term arithmetic progression. They prove that there
exists an $AP_3$ - covering sequence $A$ such that
$\limsup\limits_{n\to\infty}{A(n)}/{\sqrt n}\le  34$. In this
note, we prove that there exists an $AP_3$ - covering sequence $A$
such that $\limsup\limits_{n\to\infty}{A(n)}/{\sqrt n}=\sqrt{15}$.

 \vskip 3mm
 {\bf 2010 Mathematics Subject Classification:} 11B75

 {\bf Keywords and phrases:} Stanley sequences; $AP_3$ - covering sequences;  arithmetic
 progressions

\vskip 5mm

\section{Introduction}

Given an integer $k\ge 3$ and a set $A_0=\{ a_1,\dots , a_t\}
(a_1<\cdots <a_t)$ of nonnegative integers such that $\{ a_1,\dots
, a_t\}$ does not contain a $k$-term arithmetic progression.
Define $a_{t+1}, \dots $ by the greedy algorithm: for any $l\ge
t$, $a_{l+1}$ is the smallest integer $a
> a_{l}$ such that $\{ a_1, \dots , a_{l}, a\}$  does not contain
a $k$-term arithmetic progression. The sequence $A=\{ a_1, a_2,
\dots \} $ is called the Stanley sequence of order $k$ generated
by $A_0$. It is known that if $A$ is a Stanley sequence of order
$3$, then
$$\liminf_{n\to\infty}\frac{A(n)}{\sqrt n} \ge \sqrt 2 \quad \text{ (see
\cite{GerverRamsey} and \cite{Moy})}$$ and
$$\limsup_{n\to\infty}\frac{A(n)}{\sqrt n} \ge 1.77 \quad \text{ (see
\cite{chendai})} .$$ For related results, one may refer to
\cite{ErdosLevRauzySandor} and \cite{OdlyzkoStanley}. Recently,
Kiss, S\' andor and Yang \cite{KissSandorYang} introduced the
following notation: a sequence $A$ of nonnegative integers is
called an $AP_k$ - covering sequence if there exists an integer
$n_0$ such that if $n
> n_0$, then there exist $a_1\in A, \dots , a_{k-1}\in A$,
$a_1<a_2<\cdots <a_{k-1}<n$ such that $a_1, \dots , a_{k-1}, n$
form a $k$-term arithmetic progression. They \cite{KissSandorYang}
proved that there exists an $AP_3$ - covering sequence $A$ such
that
$$\limsup_{n\to\infty}\frac{A(n)}{\sqrt n}\le  34.$$

In this note,  the following result is proved.

\begin{theorem}\label{thm1}There exists an $AP_3$ - covering sequence $A$ such that
\begin{equation}\label{eq2}\limsup_{n\to\infty}\frac{A(n)}{\sqrt
n}=\sqrt{15}.\end{equation}
\end{theorem}

If $A$ is a Stanley sequence of order $k$, then $A$ does not
contain a $k$-term arithmetic progression. If $A$ is an $AP_k$ -
covering sequence of order $k$, then $A$  contains infinitely many
$k$-term arithmetic progressions. So none of sequences is both a
Stanley sequence of order $k$ and an $AP_k$ - covering sequence.
We pose a problem here.

\begin{problem} Is there a  Stanley sequence of order $k+1$ which is also an $AP_k$ - covering sequence?  \end{problem}

We introduce a new notation here which generalizes both Stanley
sequences of order $k$ and  $AP_k$ - covering sequences.  A
sequence $A$ of nonnegative integers is called a weak $AP_k$ -
covering sequence if there exists an integer $n_0$ such that if $n
> n_0$ and $n\notin A$, then there exist $a_1\in A, \dots , a_{k-1}\in A$,
$a_1<a_2<\cdots <a_{k-1}<n$ such that $a_1, \dots , a_{k-1}, n$
form a $k$-term arithmetic progression. Clearly,  a Stanley
sequence of order $k$ is also a weak $AP_k$ - covering sequence
and an $AP_k$ - covering sequence  of order $k$ is also a weak
$AP_k$ - covering sequence.

\section{Proof of Theorem \ref{thm1}}

Let
$$T_l=\left\{ u 4^l +\sum_{i=0}^{l-1} v_i 4^i : u\in \{ 1, 2, 3, 4
\} , v_i\in\{ 1, 2\} \right\} ,\quad l=0,1, \dots $$ and
$$A=\bigcup_{l=0}^\infty T_l.$$

First, we prove that $A$ is an $AP_3$ - covering sequence.

Let $n\ge 32$. We will prove that there exist $a, b\in A$ with
$a<b<n$ such that $a, b, n$ form a $3$-term arithmetic
progression. By $n\ge 32$,  there exists an integer $l\ge 2$ such
that $2\cdot 4^l\le n<2\cdot 4^{l+1}=8\cdot 4^l$. Let $m$ be the
integer with $m4^l\le n<(m+1)4^l$. Then $2\le m\le 7$ and
$$ 0\le n-m4^l<4^l.$$
Thus $n-m4^l$ can be written as
$$n-m4^l=\sum_{i=0}^{l-1} m_i 4^i,\quad  m_i\in \{ 0,1,2,3\} .$$
If $m_i=0$, then we take $v_{1, i} =1$ and $v_{2, i} =2$. If
$m_i\in \{ 1, 2\} $, then we take $v_{1, i} =v_{2, i} =m_i$. If
$m_i=3$, then we take $v_{1, i} =2$ and $v_{2, i} =1$.  If $m=2$,
then we take $u_1=1$ and $u_2=0$. If $m=3$, then we take $u_1=2$
and $u_2=1$. If $m=4$, then we take $u_1=2$ and $u_2=0$. If $m=5$,
then we take $u_1=3$ and $u_2=1$. If $m=6$, then we take $u_1=3$
and $u_2=0$. If $m=7$, then we take $u_1=4$ and $u_2=1$. Let
$$a=u_2 4^l +\sum_{i=0}^{l-1} v_{2, i} 4^i,\quad b=u_1 4^l +\sum_{i=0}^{l-1} v_{1, i}
4^i.$$ It is clear that $1\le a<b<n$, $a,b\in T_l\cup
T_{l-1}\subseteq A$ and $a, b, n$ form a $3$-term arithmetic
progression. Hence  $A$ is an $AP_3$ - covering sequence.

Now we prove that \eqref{eq2} holds.

Let
$$A=\{ n_1, n_2, \dots  \}, \quad n_1<n_2<\cdots .$$
For $n_j<m<n_{j+1}$, we have
$$\frac{A(m)}{\sqrt{m}}
=\frac{A(n_j)}{\sqrt{m}}<\frac{A(n_j)}{\sqrt{n_j}}.$$ It follows
that
$$\limsup_{n\to\infty}\frac{A(n)}{\sqrt n}=\limsup_{j\to\infty} \frac{A(n_j)}{\sqrt{n_j}}.$$
Let
$$n_j=u 4^l +\sum_{i=0}^{l-1} v_i 4^i,\quad u\in \{ 1, 2, 3, 4
\} , v_i\in\{ 1, 2\} \ (0\le i\le l-1) .$$ Then
\begin{equation}\label{eq1}A(n_j)=(u-1) 2^l+\sum_{i=1}^{l-1}
(v_i-1) 2^i +v_0 + 4 (2^{l-1}+\cdots +2+1).\end{equation} It is
clear that
$$n_j\ge u4^l+v_{l-1}4^{l-1}+\frac 13(4^{l-1}-1)=(4u+v_{l-1}+\frac 13)4^{l-1}-\frac 13,$$
$$A(n_j)\le (u-1)2^l+(v_{l-1}-1)2^{l-1}+2^{l-1}+4(2^l-1)=(2u+6+v_{l-1})2^{l-1}-4.$$ Since
$$2u+6+v_{l-1}<4\sqrt{4u+v_{l-1}+\frac 13} $$ for $u\in \{
1,2,3,4\}$ and $v_{l-1}\in \{ 1, 2\} $, it follows that
$$A(n_j)\le (2u+6+v_{l-1})2^{l-1}-4 < 4\sqrt{n_j+\frac
13}-4<4\sqrt{n_j}.$$ If $v_i=1$ for some $0\le i\le l-1$, then
$n_j+4^i\in A$ and by \eqref{eq1}, we have
$A(n_j+4^i)=A(n_j)+2^i$. Since $n_j>4^l\ge 4^{i+1}$, it follows
that $\sqrt{n_j+4^i}+\sqrt{n_j}>4\cdot 2^i$. That is,
$2^i>4(\sqrt{n_j+4^i}-\sqrt{n_j})$. By $A(n_j)<4\sqrt{n_j}$, we
have
$$A(n_j)(\sqrt{n_j+4^i}-\sqrt{n_j})<4\sqrt{n_j}(\sqrt{n_j+4^i}-\sqrt{n_j})<2^i \sqrt{n_j} .$$
So
$$(A(n_j)+2^i) \sqrt{n_j} >A(n_j)\sqrt{n_j+4^i} .$$
Hence
$$\frac{A(n_j+4^i)}{\sqrt{n_j+4^i}}=\frac{A(n_j)+2^i}{\sqrt{n_j+4^i}}>\frac{A(n_j)}{\sqrt{n_j}}.$$
So we need only consider those $n_j$ with all $v_i=2$. Let
$$q_{u,l}=u 4^l+\sum_{i=0}^{l-1} 2\cdot 4^i=(u+\frac 23 ) 4^l-\frac 23.$$
By \eqref{eq1}, $A(q_{u,l})=(u+4) 2^l -4$. It follows that
$$\lim_{l\to \infty } \frac{A(q_{u,l})}{\sqrt{q_{u,l}}}
=\frac{u+4}{\sqrt{u+2/3}}.$$
 Hence
$$\limsup_{n\to\infty}\frac{A(n)}{\sqrt n}=\limsup_{j\to\infty} \frac{A(n_j)}{\sqrt{n_j}}
=\max\left\{ \frac{u+4}{\sqrt{u+2/3}} : u=1,2,3,4\right\}
=\sqrt{15}.$$ This completes the proof.

\end{document}